\documentclass[11pt]{amsart}
\usepackage[latin1]{inputenc}
\usepackage[T1]{fontenc}
\usepackage{hyperref}
\usepackage{amsfonts}
\usepackage{amsmath}
\usepackage{epsf, subfigure, verbatim}
\usepackage{amssymb}
\usepackage{amsthm}
\usepackage{latexsym}
\usepackage{layout}

\newtheorem{theorem}{Theorem}

\newtheorem{corollary}{Corollary}

\newcommand{\reals}{\mathbb{R}}

\newcommand{\ind}{\mathbb{I}}
\newcommand{\e}{\mathbb{E}}

\newcommand{\ee}{\mathrm{e}}

\newcommand{\p}{\mathbb{P}}

\newcommand{\qscale}{W^{(\lambda)}}
\newcommand{\qscaleprime}{W^{(\lambda) \prime}}
\newcommand{\zscale}{Z^{(\lambda)}}

\newcommand{\up}{U^{\epsilon,b}}
\newcommand{\low}{L^{\epsilon,b}}
\newcommand{\upinf}{U^{\epsilon}}
\newcommand{\lowinf}{L^{\epsilon}}

\begin{document}

\title[Occupation times]{Occupation times of spectrally negative L\'evy processes with applications}

\author[Landriault, Renaud and Zhou]{David Landriault}
\address{Department of Statistics and Actuarial Science, University of Waterloo, 200 University Avenue West, Waterloo (Ontario) N2L 3G1, Canada}
\email{dlandria@math.uwaterloo.ca}

\author[]{Jean-Fran\c{c}ois Renaud}
\address{D\'epartement de math\'ematiques, Universit\'e du Qu\'ebec \`a Montr\'eal (UQAM), 201 av.\ Pr\'esident-Kennedy, Montr\'eal (Qu\'ebec) H2X 3Y7, Canada}
\email{renaud.jf@uqam.ca}

\author[]{Xiaowen Zhou}
\address{Department of Mathematics and Statistics, Concordia University, 1455 de Maisonneuve Blvd W., Montr\'{e}al (Qu\'{e}bec) H3G 1M8, Canada}
\email{xzhou@mathstat.concordia.ca}

\date{May 3, 2011}


\begin{abstract}
In this paper, we compute the Laplace transform of occupation times (of the negative half-line) of spectrally negative L\'evy processes. Our results are extensions of known results for standard Brownian motion and jump-diffusion processes. The results are expressed in terms of the so-called scale functions of the spectrally negative L\'evy process and its Laplace exponent. Applications to insurance risk models are also presented.
\end{abstract}

\maketitle

\section{Introduction and results}

Let $X = (X_t)_{t \geq 0}$ be a spectrally negative L\'evy process, that is a L\'{e}vy process with no positive jumps. The law of $X$ such that $X_0 = x$ is denoted by $\p_x$ and the corresponding expectation by $\e_x$. We write $\p$ and $\e$ when $x=0$. As the L\'{e}vy process $X$ has no positive jumps, its Laplace transform exists and is given by
$$
\e \left[ \mathrm{e}^{\theta X_t} \right] = \mathrm{e}^{t \psi(\theta)} ,
$$
for $\theta,t \geq 0$, where
$$
\psi(\theta) = \gamma \theta + \frac{1}{2} \sigma^2 \theta^2 + \int_{-\infty}^0 \left( \mathrm{e}^{\theta z} - 1 - \theta z \ind_{(-1,0)}(z) \right) \Pi(\mathrm{d}z) ,
$$
for $\gamma \in \reals$ and $\sigma \geq 0$. Also, $\Pi$ is a $\sigma$-finite measure on $(-\infty,0)$ such that
$$
\int_{-\infty}^0 (1 \wedge z^2) \Pi(\mathrm{d}z) < \infty .
$$
The measure $\Pi$ is called the L\'{e}vy measure of $X$, while $(\gamma,\sigma,\Pi)$ is referred to as the L\'evy triplet of $X$. Note that $\e \left[ X_1 \right] = \psi'(0+)$.

For an arbitrary spectrally negative L\'evy process, the Laplace exponent $\psi$ is strictly convex and $\lim_{\theta \to \infty} \psi(\theta) = \infty$. Thus, there exists a function $\Phi \colon [0,\infty) \to [0,\infty)$ defined by $\Phi(\theta) = \sup \{ \xi \geq 0 \mid \psi(\xi) = \theta\}$ (its right-inverse) and such that
$$
\psi ( \Phi(\theta) ) = \theta , \quad \theta \geq 0 .
$$

We first examine the total occupation time of the negative half-line $(-\infty,0)$.
%
\begin{theorem}\label{T:unrestricted}
If $\psi'(0+) > 0$, then, for $\lambda \geq 0$,
\begin{equation}\label{secondresult}
\e \left[ \mathrm{e}^{- \lambda \int_0^{\infty} \mathbb{I}_{\{X_s \leq 0\}} \mathrm{d}s } \right] = \psi'(0+) \frac{\Phi(\lambda)}{\lambda} ,
\end{equation}
where $\Phi(\lambda)/\lambda$ is to be understood in the limiting sense when $\lambda=0$.
\end{theorem}

We now recall the definition of the $q$-scale function $W^{(q)}$. For $q \geq 0$, the $q$-scale function of the process $X$ is defined as the function with Laplace transform on $[0,\infty)$ given by
$$
\int_0^{\infty} \mathrm{e}^{- \theta z} W^{(q)} (z) \mathrm{d}z = \frac{1}{\psi(\theta) - q} , \quad \text{for $\theta > \Phi(q)$,}
$$
and such that $W^{(q)}(x)=0$ for $x<0$. This function is unique, continuous, positive and strictly increasing. We write $W = W^{(0)}$ when $q=0$. We have that $W^{(q)}$ is differentiable except for at most countably many points; see \cite{loeffen2009}. Moreover, $W^{(q)}$ is continuously differentiable if $X$ has paths of unbounded variation or if the tail of the L\'evy measure, i.e., the function $x \mapsto \Pi(-\infty,-x)$ on $(-\infty,0)$, is continuous. Further, $W^{(q)}$ is twice continuously differentiable on $(0,\infty)$ if $\sigma>0$. For more details on the smoothness properties of the $q$-scale function, see \cite{chankypsavov}. We will also use the functions $\left\lbrace \overline{W}^{(q)} ; q \geq 0 \right\rbrace$ and $\left\lbrace Z^{(q)} ; q \geq 0 \right\rbrace$ defined by
$$
\overline{W}^{(q)}(x) = \int_0^x W^{(q)} (z) \mathrm{d}z .
$$
and
$$
Z^{(q)}(x) = 1 + q \overline{W}^{(q)}(x) .
$$

We can now state the following corollary to Theorem~\ref{T:unrestricted}.
\begin{corollary}\label{C:unrestricted}
If $\psi'(0+) > 0$, then, for $\lambda > 0$ and $x \geq 0$,
$$
\e_x \left[ \mathrm{e}^{- \lambda \int_0^{\infty} \mathbb{I}_{\{X_s \leq 0\}} \mathrm{d}s } \right] = \psi'(0+) \Phi(\lambda) \int_0^\infty \mathrm{e}^{-\Phi(\lambda) z} W(x+z) \mathrm{d}z ,
$$
\end{corollary}

Note that when $x=0$, by the definition of scale function, we have
$$
\int_0^\infty \mathrm{e}^{-\Phi(\lambda) z} W(z) \mathrm{d}z = \frac{1}{\psi(\Phi(\lambda)} = \frac{1}{\lambda} ,
$$
therefore recovering Theorem~\ref{T:unrestricted} as a special case.
\\

These two results generalize the work done in \cite{zhangwu2002} where the sum of a compound Poisson process and a Brownian motion is analyzed (see, e.g., Equation (4.9) in that paper). Using ruin theory terminology, they study the duration of negative surplus, also called the \textit{time in red}, in such an insurance risk model; their work is itself an extension of \cite{dosreis1993} in the pure compound Poisson case.
\\

We now examine the occupation time of $(-\infty,0)$ until a negative level $-b$ is crossed for the first time. Let $\tau_{-b}^-$ be the first passage time below $-b$ of $X$:
$$
\tau_{-b}^- = \inf \{t > 0 \colon X_t < -b \} .
$$
\begin{theorem}\label{T:restricted}
If $\psi'(0+) \geq 0$, then, for $\lambda \geq 0$,
\begin{equation}\label{mainres}
\e \left[ \mathrm{e}^{- \lambda \int_0^{\tau_{-b}^-} \mathbb{I}_{\{X_s \leq 0\}} \mathrm{d}s } \right] = \frac{\psi'(0+) + \frac{\sigma^2}{2} \frac{A_1^{(\lambda)}(b)}{\qscale(b)} + \int_{-\infty}^{0-} A_2^{(\lambda)}(x) \int_{0+}^\infty \Pi(\mathrm{d}x-y) \mathrm{d}y}{\psi'(0+) + \frac{\sigma^2}{2} \frac{\qscaleprime(b)}{\qscale(b)} + \int_{-\infty}^{0-} A_3^{(\lambda)}(x) \int_{0+}^\infty \Pi(\mathrm{d}x-y) \mathrm{d}y} ,
\end{equation}
where 
$$
A_1^{(\lambda)}(b) = \zscale(b) \qscaleprime(b) - \lambda (\qscale(b))^2,
$$
$$
A_2^{(\lambda)}(x) = \zscale(x+b) - \zscale(b) \frac{\qscale(x+b)}{\qscale(b)}
$$
and
$$
A_3^{(\lambda)}(x) = 1 - \frac{\qscale(x+b)}{\qscale(b)} .
$$
\end{theorem}
As a special case, we recover the corresponding result for standard Brownian motion; see, e.g., \cite{itomckean1974,karatzasshreve1991}.

Our proofs use fluctuation identities for spectrally negative L\'evy processes and, as a consequence, our results are expressed in terms of the so-called scale functions of the spectrally negative L\'evy process (see, e.g., \cite{bertoin1996, bertoin1997}) and its Laplace exponent. Only elementary arguments are needed. To the authors' knowledge, the literature seems rather scarce on the relationship between scale functions and certain occupation times of a general spectrally negative L\'evy process. However, similar results can be found in \cite{kadankovkadankova2004}.

The current work has been partly motivated by the study of an insurance risk model with implementation delays. Insurance risk models use stochastic processes to describe the surplus of an insurance company. In risk models of a Parisian nature, an implementation delay in the recognition of an insurer's capital insufficiency is applied. More precisely, it is assumed that ruin occurs as soon as an excursion below a critical level is longer than a deterministic time; such models have been studied only very recently in \cite{czarnapalmowski2010,dassioswu2009b,landriaultetal2010} and the idea has been borrowed from finance and Parisian barrier options (see \cite{chesneyetal1997}). Of more interest in our context is the work by the same authors: in \cite{landriaultetal2010}, instead of a deterministic delay, an exponentially distributed grace period is used in the definition of the Parisian ruin. It turns out that the probability of ruin in this model is strongly related to the occupation time of the underlying process. The reader is invited to consult Section~\ref{S:riskmodels} to obtain further details on this connection.

The rest of the paper is organized as follows. In the next Section, we recall the relevant notions and results on scale functions and fluctuation identities. Then, in Section 3 and 4, we prove Theorem~\ref{T:restricted} and Theorem~\ref{T:unrestricted} respectively. Section 5 presents applications to the case of a Brownian motion with drift and insurance risk models.

\section{Scale functions and fluctuation identities}\label{S:preliminaries}

We recall some of the properties of the $q$-scale function $W^{(q)}$ and its use in fluctuation theory. Let $d=\gamma-\int_{-1}^0 z \Pi(\mathrm{d}z)$. The initial values of $W^{(q)}$ and $W^{(q)\prime}$ are known to be
\begin{equation*}
W^{(q)}(0+)=
\begin{cases}
1/d & \text{when $\sigma=0$ and $\int_{-1}^0 z \Pi(\mathrm{d}z) < \infty$},  \\
0 & \text{otherwise},
\end{cases}
\end{equation*}
and
\begin{equation*}
W^{(q)\prime}(0+)=
\begin{cases}
2/\sigma^2 & \text{when $\sigma>0$,} \\
(\Pi(-\infty,0)+q)/d^2 & \text{when $\sigma=0$ and $\Pi(-\infty,0)<\infty$,} \\
\infty & \text{otherwise.}
\end{cases}
\end{equation*}

Now, define
$$
\tau_0^- = \inf \{t > 0 \colon X_t < 0 \},
$$
and, for $a > 0$,
$$
\tau_a^+ = \inf \{t > 0 \colon X_t > a \} ,
$$
with the convention $\inf\emptyset=\infty$. It is well known (see, e.g., \cite{kyprianou2006}) that, for $x \leq a$,
$$
\e_x \left[ \ee^{-q \tau_a^+} ; \tau_a^+ < \tau_0^- \right] = \frac{W^{(q)}(x)}{W^{(q)}(a)} ,
$$
$$
\e_x \left[ \ee^{-q \tau_0^-} ; \tau_0^- < \tau_a^+ \right] = Z^{(q)}(x) - Z^{(q)}(a) \frac{W^{(q)}(x)}{W^{(q)}(a)} ,
$$
and
$$
\e_x \left[ \ee^{-q \tau_a^+} ; \tau_a^+ < \infty \right] = \ee^{-\Phi(q) (a-x)} .
$$
If $\psi'(0+) \geq 0$, then $\p_x \left\lbrace \tau_a^+ < \infty \right\rbrace = 1$ and therefore
$$
\e_x \left[ \ee^{-q \tau_a^+} \right] = \ee^{-\Phi(q) (a-x)} .
$$
Also, we have that
\begin{multline}\label{E:deficit}
\p_x \left\lbrace \tau_0^- < \infty , X_{\tau_0^-} \in \mathrm{d}z \right\rbrace = \frac{\sigma^2}{2} \left[ W'(x) - \Phi(0) W(x) \right] \delta_0( \mathrm{d}z) \\
+ \int_{0+}^\infty \Pi(\mathrm{d}z-y) \left\lbrace \ee^{-\Phi(0) y} W(x) - W(x-y) \right\rbrace \mathrm{d}y ,
\end{multline}
where $\delta_0$ is the Dirac measure at $0$. The first term of this measure corresponds to the case when $X_{\tau_0^-} = 0$, a behaviour called \textit{creeping}. Note that we write $\int_{0+}^\infty$ to mean $\int_{(0,\infty)}$; similarly, we will write $\int_{-\infty}^{0-}$ to mean $\int_{(-\infty,0)}$, $\int_{-\infty}^{0+}$ to mean $\int_{(-\infty,0]}$, etc.\\

Using the distribution in Equation~\eqref{E:deficit}, together with the fact that $\p_x \left\lbrace \tau_0^- < \infty \right\rbrace = 1 - \psi'(0+) W(x)$ if $\psi'(0+) \geq 0$ (in which case $\Phi(0)=0$), we obtain
\begin{equation}\label{E:trick_ruinproba}
1 = \psi'(0+)W(x) + \frac{\sigma^2}{2} W'(x) + \int_{-\infty}^{0-} \int_{0+}^\infty \Pi(\mathrm{d}z-y) \left\lbrace W(x) - W(x-y) \right\rbrace \mathrm{d}y .
\end{equation}

For more details on spectrally negative L\'{e}vy processes and fluctuation identities, the reader is referred to \cite{kyprianou2006}. Further information, examples and numerical techniques related to the computation of scale functions can be found in \cite{chankypsavov,egamiyamazaki2010,hubalekkyprianou2010,kyprianourivero2008,surya2008}.

%
%
%

%
\section{Proof of Theorem~\ref{T:restricted}}

The main idea of the proof consists in defining a quantity underestimating and overestimating the occupation time
$$
\int_0^{\tau_{-b}^-} \mathbb{I}_{\{X_s \leq 0\}} \mathrm{d}s .
$$
This respectively leads to an upper and a lower bound to its Laplace transform. Subsequently, by taking an appropriate limit, we show that the two bounds converge to the expression on the right-hand side of \eqref{mainres}. 

First, we provide a lower bound to this Laplace transform by overestimating the occupation time. To this end, we consider a clock which starts at time $0$ and stops when level $\epsilon$ is attained or when level $-b$ is crossed. Then, if level $\epsilon$ was attained first, every time we go below $0$, we re-start the clock and subsequently stop it when we get back to $\epsilon$ (without going below $-b$); let $\low$ be the Laplace transform of this overestimating quantity of the occupation time when the process $X$ sits at level $\epsilon$ at time $0$. Hence, by the strong Markov property of $X$, we have 
$$
\e \left[ \mathrm{e}^{- \lambda \int_0^{\tau_{-b}^-} \mathbb{I}_{\{X_s \leq 0\}} \mathrm{d}s } \right] \geq \e \left[ \mathrm{e}^{- \lambda \tau_\epsilon^+} ; \tau_\epsilon^+ < \tau_{-b}^- \right] \low = \frac{\qscale(b)}{\qscale(b+\epsilon)} \low .
$$
Using the strong Markov property and the spatial homogeneity of $X$, we get
\begin{multline*}
\low = \p_\epsilon \{ \tau_0^- = \infty \} \\
+ \int_{-\infty}^{0+} \p_\epsilon \left\lbrace \tau_0^- < \infty , X_{\tau_0^-} \in \mathrm{d}x \right\rbrace \left\lbrace \e_x \left[ \mathrm{e}^{- \lambda \tau_{-b}^-} ; \tau_{-b}^- < \tau_\epsilon^+ \right] + \e_x \left[ \mathrm{e}^{- \lambda \tau_\epsilon^+} ; \tau_\epsilon^+ < \tau_{-b}^- \right] \low \right\rbrace .
\end{multline*}
Note that when $x < -b$,
$$
\e_x \left[ \mathrm{e}^{- \lambda \tau_{-b}^-} ; \tau_{-b}^- < \tau_\epsilon^+ \right] = 1
$$
and
$$
\e_x \left[ \mathrm{e}^{- \lambda \tau_\epsilon^+} ; \tau_\epsilon^+ < \tau_{-b}^- \right] = 0 .
$$
Consequently,
$$
\low = \frac{\p_\epsilon \{ \tau_0^- = \infty \} + \int_{-\infty}^{0+} \p_\epsilon \left\lbrace \tau_0^- < \infty , X_{\tau_0^-} \in \mathrm{d}x \right\rbrace \e_x \left[ \mathrm{e}^{- \lambda \tau_{-b}^-} ; \tau_{-b}^- < \tau_\epsilon^+ \right]}{1 - \int_{-\infty}^{0+} \p_\epsilon \left\lbrace \tau_0^- < \infty , X_{\tau_0^-} \in \mathrm{d}x \right\rbrace \e_x \left[ \mathrm{e}^{- \lambda \tau_\epsilon^+} ; \tau_\epsilon^+ < \tau_{-b}^- \right]} = \frac{\low_1 + \low_2}{\low_3} ,
$$
where, using some of the fluctuation identities in Section~\ref{S:preliminaries},
$$
\low_1 = \psi'(0+) W(\epsilon) + \frac{\sigma^2}{2} W'(\epsilon) \left[ \zscale(b) - \zscale(b+\epsilon) \frac{\qscale(b)}{\qscale(b+\epsilon)} \right] ,
$$
$$
\low_2 = \int_{-\infty}^{0-} \left[ \zscale(x+b) - \zscale(b+\epsilon) \frac{\qscale(x+b)}{\qscale(b+\epsilon)} \right] \int_{0+}^\infty \Pi(\mathrm{d}x-y) \left\lbrace W(\epsilon) - W(\epsilon-y) \right\rbrace \mathrm{d}y ,
$$
and,
$$
\low_3 = 1 - \frac{\sigma^2}{2} W'(\epsilon) \frac{\qscale(b)}{\qscale(b+\epsilon)} - \int_{-\infty}^{0-} \frac{\qscale(x+b)}{\qscale(b+\epsilon)} \int_{0+}^\infty \Pi(\mathrm{d}x-y) \left\lbrace W(\epsilon) - W(\epsilon-y) \right\rbrace \mathrm{d}y.
$$
Using Equation~\eqref{E:trick_ruinproba}, we get
\begin{multline*}
\low_3 = \psi'(0+) W(\epsilon) + \frac{\sigma^2}{2} \frac{W'(\epsilon)}{\qscale(b+\epsilon)} \left[ \qscale(b+\epsilon)-\qscale(b) \right] \\
+ \int_{-\infty}^{0-} \left[ 1 - \frac{\qscale(x+b)}{\qscale(b+\epsilon)} \right] \int_{0+}^\infty \Pi(\mathrm{d}x-y) \left\lbrace W(\epsilon) - W(\epsilon-y) \right\rbrace \mathrm{d}y .
\end{multline*}

Alternatively, we now develop a scheme to underestimate the occupation time in question. Every time we go below $-\epsilon$ ($-\epsilon > -b$), we start the clock and stop it when we get back to $0$ (without going below $-b$); let $\up$ be the Laplace transform of this underestimating quantity of the occupation time when $X_{0}=0$. Hence, by the strong Markov property of $X$, we also have 
$$
\e \left[ \mathrm{e}^{- \lambda \int_0^{\tau_{-b}^-} \mathbb{I}_{\{X_s \leq 0\}} \mathrm{d}s } \right] \leq \up .
$$
As above, using the strong Markov property and the spatial homogeneity of $X$, we can write
\begin{multline*}
\up = \p_{\epsilon} \{ \tau_0^- = \infty \} \\
+ \int_{-\infty}^{0+} \p_\epsilon \left\lbrace \tau_0^- < \infty , X_{\tau_0^-} \in \mathrm{d}x \right\rbrace \left\lbrace \e_x \left[ \mathrm{e}^{- \lambda \tau_{-b+\epsilon}^-} ; \tau_{-b+\epsilon}^- < \tau_\epsilon^+ \right] + \e_x \left[ \mathrm{e}^{- \lambda \tau_\epsilon^+} ; \tau_\epsilon^+ < \tau_{-b+\epsilon}^- \right] \up \right\rbrace
\end{multline*}
and then
$$
\up = \frac{\up_1 + \up_2}{\up_3} ,
$$
where
$$
\up_1 = \psi'(0+) W(\epsilon) + \frac{\sigma^2}{2} W'(\epsilon) \left[ \zscale(b-\epsilon) - \zscale(b) \frac{\qscale(b-\epsilon)}{\qscale(b)} \right] ,
$$
$$
\up_2 = \int_{-\infty}^{0-} \left[ \zscale(x+b-\epsilon) - \zscale(b) \frac{\qscale(x+b-\epsilon)}{\qscale(b)} \right] \int_{0+}^\infty \Pi(\mathrm{d}x-y) \left\lbrace W(\epsilon) - W(\epsilon-y) \right\rbrace \mathrm{d}y ,
$$
and
$$
\up_3 = 1 - \frac{\sigma^2}{2} W'(\epsilon) \frac{\qscale(b-\epsilon)}{\qscale(b)} - \int_{-\infty}^{0-} \frac{\qscale(x+b-\epsilon)}{\qscale(b)} \int_{0+}^\infty \Pi(\mathrm{d}x-y) \left\lbrace W(\epsilon) - W(\epsilon-y) \right\rbrace \mathrm{d}y .
$$
As for $\low_3$, using Equation~\eqref{E:trick_ruinproba}, we get
\begin{multline*}
\up_3 = \psi'(0+) W(\epsilon) + \frac{\sigma^2}{2} \frac{W'(\epsilon)}{\qscale(b)} \left[ \qscale(b)-\qscale(b-\epsilon) \right] \\
+ \int_{-\infty}^{0-} \left[ 1 - \frac{\qscale(x+b-\epsilon)}{\qscale(b)} \right] \int_{0+}^\infty \Pi(\mathrm{d}x-y) \left\lbrace W(\epsilon) - W(\epsilon-y) \right\rbrace \mathrm{d}y .
\end{multline*}

%
\subsection{Proof if $\sigma > 0$}

First, we assume that $X$ has a Brownian component, that is $\sigma > 0$. In this case, $W(0)=0$ and $\qscaleprime (0+) < \infty$. Then,
\begin{align*}
\lim_{\epsilon \to 0} \frac{\low_1}{\epsilon} &= \lim_{\epsilon \to 0} \psi'(0+) \frac{W(\epsilon)}{\epsilon} + \frac{\sigma^2}{2} \frac{W'(\epsilon)}{\qscale(b+\epsilon)} \left[ \frac{\zscale(b) \qscale(b+\epsilon) - \zscale(b+\epsilon) \qscale(b)}{\epsilon} \right] \\
&= \psi'(0+) W'(0+) + \frac{\sigma^2}{2} \frac{W'(0+)}{\qscale(b)} \left[ \zscale(b) \qscaleprime(b) - \lambda (\qscale(b))^2 \right] ,
\end{align*}
\begin{align*}
\lim_{\epsilon \to 0} \frac{\low_2}{\epsilon} &= \lim_{\epsilon \to 0} \int_{-\infty}^{0-} \left[ \zscale(x+b) - \zscale(b+\epsilon) \frac{\qscale(x+b)}{\qscale(b+\epsilon)} \right] \int_{0+}^\infty \Pi(\mathrm{d}x-y) \left\lbrace \frac{W(\epsilon)}{\epsilon} - \frac{W(\epsilon-y)}{\epsilon} \right\rbrace \mathrm{d}y \\
&= W'(0+) \int_{-\infty}^{0-} \left[ \zscale(x+b) - \zscale(b) \frac{\qscale(x+b)}{\qscale(b)} \right] \int_{0+}^\infty \Pi(\mathrm{d}x-y) \mathrm{d}y
\end{align*}
and, similarly,
\begin{align*}
\lim_{\epsilon \to 0} \frac{\low_3}{\epsilon} &= \psi'(0+) \lim_{\epsilon \to 0} \frac{W(\epsilon)}{\epsilon} + \frac{\sigma^2}{2} \lim_{\epsilon \to 0} \frac{W'(\epsilon)}{\qscale(b+\epsilon)} \left[ \frac{\qscale(b+\epsilon)-\qscale(b)}{\epsilon} \right] \\
& \qquad + \lim_{\epsilon \to 0} \int_{-\infty}^{0-} \left[ 1 - \frac{\qscale(x+b)}{\qscale(b+\epsilon)} \right] \int_{0+}^\infty \Pi(\mathrm{d}x-y) \left\lbrace \frac{W(\epsilon)}{\epsilon} - \frac{W(\epsilon-y)}{\epsilon} \right\rbrace \mathrm{d}y \\
&= \psi'(0+) W'(0+) + \frac{\sigma^2}{2} W'(0+) \frac{\qscaleprime(b)}{\qscale(b)} + W'(0+) \int_{-\infty}^{0-} \left[ 1 - \frac{\qscale(x+b)}{\qscale(b)} \right] \int_{0+}^\infty \Pi(\mathrm{d}x-y) \mathrm{d}y .
\end{align*}
The result follows.

\subsection{Proof if $\sigma = 0$}

If $X$ does not have a Brownian component, we adapt the previous method as follows:
$$
\frac{\low_1}{W(\epsilon)} = \psi'(0+) , \, \text{for all $\epsilon > 0$,}
$$
\begin{align*}
\lim_{\epsilon \to 0} \frac{\low_2}{W(\epsilon)} &= \lim_{\epsilon \to 0} \int_{-\infty}^{0-} \left[ \zscale(x+b) - \zscale(b+\epsilon) \frac{\qscale(x+b)}{\qscale(b+\epsilon)} \right] \int_{0+}^\infty \Pi(\mathrm{d}x-y) \left\lbrace \frac{W(\epsilon)}{W(\epsilon)} - \frac{W(\epsilon-y)}{W(\epsilon)} \right\rbrace \mathrm{d}y \\
&= \lim_{\epsilon \to 0} \int_{-\infty}^{0-} \left[ \zscale(x+b) - \zscale(b+\epsilon) \frac{\qscale(x+b)}{\qscale(b+\epsilon)} \right] \int_{0+}^\infty \Pi(\mathrm{d}x-y) \mathrm{d}y \\
&= \int_{-\infty}^{0-} \left[ \zscale(x+b) - \zscale(b) \frac{\qscale(x+b)}{\qscale(b)} \right] \int_{0+}^\infty \Pi(\mathrm{d}x-y) \mathrm{d}y
\end{align*}
and, similarly,
\begin{align*}
\lim_{\epsilon \to 0} \frac{\low_3}{W(\epsilon)} &= \psi'(0+) + \lim_{\epsilon \to 0} \int_{-\infty}^{0-} \left[ 1 - \frac{\qscale(x+b)}{\qscale(b+\epsilon)} \right] \int_{0+}^\infty \Pi(\mathrm{d}x-y) \left\lbrace \frac{W(\epsilon)}{W(\epsilon)} - \frac{W(\epsilon-y)}{W(\epsilon)} \right\rbrace \mathrm{d}y \\
&= \psi'(0+) + \int_{-\infty}^{0-} \left[ 1 - \frac{\qscale(x+b)}{\qscale(b)} \right] \int_{0+}^\infty \Pi(\mathrm{d}x-y) \mathrm{d}y .
\end{align*}
The result follows.

\subsection{Conclusion of the proof}

In all cases, i.e., when $X$ has or does not have a Brownian component, the limiting results for $\up$ can be obtained in a similar fashion. In fact, it can easily be proved that, for each $i=1,2,3$,
$$
\lim_{\epsilon \to 0} \frac{\up_i}{\low_i} = 1
$$
and then
$$
\lim_{\epsilon \to 0} \frac{\qscale(b)}{\qscale(b+\epsilon)} \low = \e \left[ \mathrm{e}^{- \lambda \int_0^{\tau_{-b}^-} \mathbb{I}_{\{X_s \leq 0\}} \mathrm{d}s } \right] = \lim_{\epsilon \to 0} \up .
$$
The details are left to the reader.

%
\section{Proof of Theorem~\ref{T:unrestricted}}

We use the same general idea as in the proof of Theorem~\ref{T:restricted}. Define $\lowinf = L^{\epsilon,\infty}$. In this case, using the strong Markov property twice, we obtain
\begin{align*}
\lowinf &= \p_\epsilon \{ \tau_0^- = \infty \} + \lowinf \e \left[ \mathrm{e}^{- \lambda \tau_\epsilon^+} \right] \int_{-\infty}^{0+} \p_\epsilon \left\lbrace \tau_0^- < \infty , X_{\tau_0^-} \in \mathrm{d}x \right\rbrace \e_x \left[ \mathrm{e}^{- \lambda \tau_0^+} \right] \\
&= \psi'(0+) W(\epsilon) + \lowinf \mathrm{e}^{-\Phi(\lambda) \epsilon} \int_{-\infty}^{0+} \p_\epsilon \left\lbrace \tau_0^- < \infty, X_{\tau_0^-} \in \mathrm{d}x \right\rbrace \mathrm{e}^{\Phi(\lambda) x} .
\end{align*}
Since, for $r > 0$,
\begin{equation}\label{E:laplace_deficit}
\e_x \left[ \mathrm{e}^{r X_{\tau_0^-} } ; \tau_0^- < \infty \right] = \mathrm{e}^{r x} - \psi(r) \mathrm{e}^{r x} \int_0^x \mathrm{e}^{-r z} W(z) \mathrm{d}z - \frac{\psi(r)}{r} W(x) ,
\end{equation}
we can write
\begin{align*}
\int_{-\infty}^{0+} \p_\epsilon \left\lbrace \tau_0^- < \infty , X_{\tau_0^-} \in \mathrm{d}x \right\rbrace \mathrm{e}^{\Phi(\lambda) x} &= \e_\epsilon \left[ \mathrm{e}^{\Phi(\lambda) X_{\tau_0^-}} ; \tau_0^- < \infty \right] \\
&= \mathrm{e}^{\Phi(\lambda) \epsilon} - \lambda \mathrm{e}^{\Phi(\lambda) \epsilon} \int_0^\epsilon \mathrm{e}^{-\Phi(\lambda) x} W(x) \mathrm{d}x - \frac{\lambda}{\Phi(\lambda)} W(\epsilon).
\end{align*}
It is easily shown that
$$
\lim_{\epsilon \to 0} \frac{\int_0^\epsilon \mathrm{e}^{-\Phi(\lambda) x} W(x) \mathrm{d}x}{W(\epsilon)} = 0.
$$
Indeed, if $W(0+)>0$ this is trivial and if $W(0+)=0$ one can use L'H\^opital's rule. Consequently,
$$
\lim_{\epsilon \to 0} \lowinf = \lim_{\epsilon \to 0} \frac{\psi'(0+) W(\epsilon)}{1- \left[ 1 - \lambda \int_0^\epsilon \mathrm{e}^{-\Phi(\lambda) x} W(x) \mathrm{d}x - \frac{\lambda}{\Phi(\lambda)} \mathrm{e}^{-\Phi(\lambda) \epsilon} W(\epsilon) \right] } = \psi'(0+) \frac{\Phi(\lambda)}{\lambda} .
$$
Similarly, if we define $\upinf = U^{\epsilon,\infty}$, one can show that
$$
\lim_{\epsilon \to 0} \upinf = \psi'(0+) \frac{\Phi(\lambda)}{\lambda}
$$
and the result immediately follows.

\section{Proof of Corollary~\ref{C:unrestricted}}

Using the strong Markov property of $X$ twice, using some of the fluctuation identities in Section~\ref{S:preliminaries} and then using Theorem~\ref{T:unrestricted}, we have
\begin{align*}
\e_x \left[ \mathrm{e}^{- \lambda \int_0^{\infty} \mathbb{I}_{\{X_s \leq 0\}} \mathrm{d}s } \right] &= \p_x \left\lbrace \tau_0^- = \infty \right\rbrace + \e_x \left[ \mathrm{e}^{- \lambda \int_0^{\infty} \mathbb{I}_{\{X_s \leq 0\}} \mathrm{d}s } ; \tau_0^- < \infty\right] \\
&= \psi'(0+) W(x) + \psi'(0+) \frac{\Phi(\lambda)}{\lambda} \int_{-\infty}^{0+} \p_x \left\lbrace \tau_0^- < \infty , X_{\tau_0^-} \in \mathrm{d}z \right\rbrace \e_z \left[ \mathrm{e}^{- \lambda \tau_0^+} ; \tau_0^+ < \infty \right] \\
&= \psi'(0+) W(x) + \psi'(0+) \frac{\Phi(\lambda)}{\lambda} \int_{-\infty}^{0+} \p_x \left\lbrace \tau_0^- < \infty , X_{\tau_0^-} \in \mathrm{d}z \right\rbrace \mathrm{e}^{\Phi(\lambda) z} \\
&= \psi'(0+) W(x) + \psi'(0+) \frac{\Phi(\lambda)}{\lambda} \e_x \left[ \mathrm{e}^{\Phi(\lambda) X_{\tau_0^-} } ; \tau_0^- < \infty \right] .
\end{align*}
Using the result in Equation~\eqref{E:laplace_deficit} once more, we get
\begin{align*}
\e_x \left[ \mathrm{e}^{- \lambda \int_0^{\infty} \mathbb{I}_{\{X_s \leq 0\}} \mathrm{d}s } \right] &= \psi'(0+) \frac{\Phi(\lambda)}{\lambda} \left\lbrace \mathrm{e}^{\Phi(\lambda) x} \left(  1 - \lambda \int_0^x \mathrm{e}^{-\Phi(\lambda) z} W(z) \mathrm{d}z \right) \right\rbrace \\
&= \psi'(0+) \frac{\Phi(\lambda)}{\lambda} \left\lbrace \lambda \int_0^\infty \mathrm{e}^{-\Phi(\lambda) z} W(x+z) \mathrm{d}z \right\rbrace ,
\end{align*}
where in the last step a change of variables and an integration by parts were undertaken. The result follows.

\section{Applications}

\subsection{Brownian motion with drift}

For Brownian motion with drift, i.e., if $X$ is of the form $X_t = mt + \sigma B_t$ with $m=\psi'(0+) \geq 0$ and $\sigma > 0$, where $B$ is a standard Brownian motion, then $\Pi \equiv 0$ and, using Theorem~\ref{T:restricted}, we have
$$
\e \left[ \mathrm{e}^{- \lambda \int_0^{\tau_{-b}^-} \mathbb{I}_{\{X_s \leq 0\}} \mathrm{d}s } \right] = \frac{m \qscale(b) + (\sigma^2/2) \left[ \zscale(b) \qscaleprime(b) - \lambda (\qscale(b))^2 \right]}{m \qscale(b) + (\sigma^2/2) \qscaleprime(b)} ,
$$
where, as one can easily verify from the definition,
$$
\qscale (x) = \frac{2}{\sqrt{m^2 + 2{\lambda} \sigma^2}} \mathrm{e}^{-(m/\sigma^2)x} \sinh \left( (x/\sigma^2) \sqrt{m^2 + 2{\lambda} \sigma^2} \right) ,
$$
for $x \geq 0$.

Then, for standard Brownian motion, in which case $m=\psi'(0+)=0$ and $\sigma=1$, we have, from Theorem~\ref{T:restricted}, that
$$
\e \left[ \mathrm{e}^{- \lambda \int_0^{\tau_{-b}^-} \mathbb{I}_{\{X_s \leq 0\}} \mathrm{d}s } \right] = \frac{\zscale(b) \qscaleprime(b) - \lambda (\qscale(b))^2}{\qscaleprime(b)} ,
$$
%
%
where
$$
\qscale (x) = \sqrt{\frac{2}{\lambda}} \sinh (x \sqrt{2 \lambda} )
$$
and then
$$
\zscale(x) = \cosh (x \sqrt{2 \lambda} ) .
$$
Recalling that $\cosh^2(x)-\sinh^2(x)=1$, we get that
\begin{align*}
\e \left[ \mathrm{e}^{- \lambda \int_0^{\tau_{-b}^-} \mathbb{I}_{\{X_s \leq 0\}} \mathrm{d}s } \right] &= \frac{2 \cosh^2 (b \sqrt{2 \lambda}) - \lambda (2/\lambda) \sinh^2 (b \sqrt{2 \lambda} ) }{2 \cosh (b \sqrt{2 \lambda} ) } \\
&= \frac{1}{\cosh (b \sqrt{2 \lambda} ) } ,
\end{align*}
therefore recovering Proposition 4.12 of \cite{karatzasshreve1991} (see also \cite{itomckean1974}).

\subsection{Insurance risk models}\label{S:riskmodels}

Classical insurance risk models describe the surplus process of an insurance company using a compound Poisson process or a Brownian motion with drift, i.e., special cases of spectrally negative L\'evy processes. In those models, it is usually assumed that the \textit{net profit condition} holds; this condition ensures that ruin will not occur almost surely. For a L\'{e}vy insurance risk process, i.e., a spectrally negative L\'evy processw with non-monotone paths, it amounts to $\e [X_1] = \psi'(0+) > 0$ and it is also equivalent to $\lim_{t \to \infty} X_t = \infty$ almost surely. An interpretation of L\'{e}vy insurance risk models for the surplus modelling of large insurance companies is for instance given in \cite{kluppelbergkyprianou2006}.

Classically, the probability of ruin has been the most studied risk measure to evaluate the quality of a company. More recently, the analysis of the duration of the negative surplus, in other words the occupation time of the negative half-line, has also been considered in the compound Poisson case \cite{dosreis1993} and then in a jump-diffusion model \cite{zhangwu2002}; see also \cite{biardetal2010}. Recall that Theorem~\ref{T:unrestricted} extends Equation (4.9) in \cite{zhangwu2002} where the sum of a compound Poisson process and a Brownian motion is considered.

%

%

We now want to provide another link between Theorem~\ref{T:unrestricted} and insurance risk models. In \cite{landriaultetal2010}, a \textit{new} definition of the time to ruin is proposed. In that paper, each excursion of the surplus process $X$ below $0$ is accompanied by an independent copy of an independent (of $X$) and exponentially distributed random variable $\ee_d$ with mean $1/d$; we will refer to it as the implementation clock. If the duration of a given excursion below $0$ is less than its associated implementation clock, then ruin does not occur. More precisely, we assume that ruin occurs at the first time $\tau_d$ that an implementation clock rings before the end of its corresponding excursion below $0$. It is worth pointing out that the time to ruin $\tau_d$ is easily defined when the L\'evy insurance risk processes $X$ has sample paths of bounded variation.

Therefore, in \cite{landriaultetal2010}, the case of a surplus process of bounded variation is considered. In that model, one can show that the probability of ruin in this Parisian risk model with exponential implementation delays can be expressed as follows:
$$
\p \left\lbrace \tau_d < \infty \right\rbrace = 1 - \e \left[ \mathrm{e}^{- d \int_0^{\infty} \mathbb{I}_{\{X_s \leq 0\}} \mathrm{d}s } \right] ,
$$
if the exponential clock has mean $1/d$. Indeed, for each excursion of length $T$, the probability to survive, i.e., the probability that the exponential random variable associated with it is larger than $T$, is equal to $\exp \{- d \, T \}$. Using the independence assumption between the clocks and summing up over all the excursions (there are countably many of them), we get
$$
\e \left[ \mathrm{e}^{- d \int_0^{\infty} \mathbb{I}_{\{X_s \leq 0\}} \mathrm{d}s } \right] .
$$
Using Theorem~\ref{T:unrestricted}, we recover the corresponding expression in \cite{landriaultetal2010}, that is
$$
\p \left\lbrace \tau_d < \infty \right\rbrace = 1 - \psi'(0+) \frac{\Phi(d)}{d} ,
$$
when the net profit condition $\psi'(0+)>0$ is verified.

More generally, using It\^o's excursion theory for spectrally negative L\'evy processes, this Parisian risk model with exponential implementation delays can also be defined when the underlying surplus process has paths of unbounded variation. It suffices to \textit{mark} the Poisson point process of excursions away from zero (see \cite{oblojpistorius2009} for a definition of the corresponding excursion process) with independent copies of the generic random variable $\ee_d$, similar to the proof of Theorem 6.16 in \cite{kyprianou2006}; for an excursion away from zero starting above zero, this time spent above zero is simply ignored. As a consequence, Theorem~\ref{T:unrestricted} provides a generalization of the probability of ruin in a general L\'evy insurance risk model with exponential implementation delays. 

\section{Acknowledgements}

Funding in support of this work was provided by the Natural Sciences and Engineering Research Council of Canada (NSERC) and the Insti\-tut de finance math\'ematique de Montr\'eal (IFM2).

%
%
\bibliographystyle{abbrv}
\bibliography{../occupation.bib}

\end{document}